\documentclass [11pt] {amsart}
\usepackage{amsthm,amsmath,amssymb,amscd,graphics,enumerate}
\usepackage{latexsym}
\usepackage{verbatim}
\usepackage{enumerate}
\usepackage{amsthm}
\usepackage{amsmath}
\usepackage{amscd}

\newcounter{theorem}[section]
\numberwithin{equation}{section}

\newtheorem{Thm}[theorem]{Theorem}
{\theoremstyle{remark}
 
\newtheorem{Rem}[theorem]{\text{\textbf{Remark}}} }
\newtheorem{Def}[theorem]{Definition}
\newtheorem{Lem}[theorem]{Lemma}
\newtheorem{Prop}[theorem]{Proposition}

{\theoremstyle{definition}
 }

\newcommand{\CC}{{\mathbb{C}}}
\newcommand{\QQ}{{\mathbb{Q}}}

\newcommand{\mO}{\mathcal{O}}
\newcommand{\mX}{\mathcal{X}}
\newcommand{\mY}{\mathcal{Y}}
\newcommand{\mW}{\mathcal{W}}

\newcommand{\mE}{\mathcal{E}}

\newcommand{\mV}{\mathcal{V}}
\newcommand{\mZ}{\mathcal{Z}}

\title{A note on derived McKay correspondence}
\author{Jiun-Cheng Chen}
\address{Department of Mathematics\\ Northwestern University\\ 2033 Sheridan Road\\ Evanston\\ IL 60208-2370\\ USA}
\curraddr{Department of Mathematics \\National Tsing Hua University\\3rd General Building, 101 Section 2, Kuang Fu Road\\Hsinchu 300\\ Taiwan}

\email{jcchen@math.northwestern.edu, jcchen@math.nthu.edu.tw}
\thanks{Research of J.-C. C. supported in part by National Center of Theoretical Science in Hsinchu, Taiwan and National Tsing Hua University Center for Mathematical Science.}

\author{Hsian-Hua Tseng}
\address{Department of Mathematics\\ University of British Columbia\\ 1984 Mathematics Road\\Vancouver\\ BC V6T 1Z1\\ Canada}
\email{hhtseng@math.ubc.ca}
\thanks{Research of H.-H. T. supported in part by Clay Mathematics Institute Liftoff Fellowship program and a postdoctoral fellowship at the Mathematical Sciences Research Institute.}
\date{\today}

\begin{document}
\begin{abstract}
We obtain a global version and a twisted version (in the sense of \cite{bp05}) of the main theorem of \cite{bkr}.
\end{abstract}
\maketitle
\section{Introduction}
We work over the field of complex numbers.

Let $X$ be an irreducible projective variety of dimension $n$. Assume that $X$ has only quotient singularities. According to \cite{vi89}, there is a smooth Deligne-Mumford stack $\mX$ with coarse moduli space $X$, such that $\mX$ and $X$ are isomorphic in codimension one. Note that $\mX$ is a quotient stack. Let $\pi: \mX\to X$ denote the projection.

In the case where $X=M/G$ with $G$ a finite group, let $Y\subset \text{G-Hilb}(M)$ be the irreducible component of the $G$-Hilbert scheme of $M$ that contains the free orbits. There is a morphism $Y\to X$ called the Hilbert-Chow morphism. The main result of \cite{bkr} can be stated as follows: Suppose that $dim Y \times_X Y\leq n+1$, then $Y$ is smooth and there is an equivalence of derived categories $D^b(Y)\simeq D^b(\mX)$.

We study the global version of this problem (that is, $X$ is not necessarily of the form $M/G$ as above). Recall that for any Deligne-Mumford stack $\mX$, \'etale locally on its coarse moduli space, $\mX$ is of the form $[U/G]$ with $U$ a scheme and $G$ a finite group. It is tempting to obtain a global crepant resolution by patching the local ones. This is, however, not obvious at all:   Suppose
 that $\{V_i\}$ is an \'etale cover of $X$, $V_i\times_X \mX \simeq [U_i/G_i]$, and for each $i$ there is a crepant resolution $\phi_i: Y_i\to V_i$. Then it is not clear that we can patch these $\{Y_i\}$ together, since
crepant resolution in higher dimension is not unique.
In dimension 3, a global crepant resolution can be built from local ones since flops preserve smoothness and two crepant resolutions can be connected by a sequence of flops, see Proposition \ref{3dcr} below. In  higher dimensions, this argument doesn't work since flops do not preserve smoothness  and may not
terminate.

An observation, which we learned from D. Abramovich, is that a certain Hilbert functor studied by Olsson-Starr \cite{os} is a good replacement for $G$-Hilbert
 schemes in the global situation. We denote the scheme representing
this Hilbert functor by $\text{Hilb}(\mX)$. We will only be interested in a particular component $\text{Hilb}'(\mX) \subset \text{Hilb}(\mX)$: Let $U \subset \mX$ be the open set of non-stacky points in $\mX$. There is a natural inclusion   $U \subset \text{Hilb}(\mX)$. The scheme $\text{Hilb}'(\mX)$ is the component which contains $U$.
There is a morphism $$\text{Hilb}'(\mX)\to X$$ induced by the functor $\pi_*$.
\begin{Thm}\label{gbkr}
Assume that $dim \text{Hilb}'(\mX)\times_X \text{Hilb}'(\mX) \leq n+1$, then the normalization $(\text{Hilb}'(\mX))^n$ is smooth and the Fourier-Mukai type transformation $$F: D^b((\text{Hilb}'(\mX))^n)\to D^b(\mX)$$ induced by the structure sheaf of the  universal object over $(\text{Hilb}'(\mX))^n$ is an equivalence of derived categories.
\end{Thm}

\begin{Rem}
Crepant resolutions, if exist, may not be unique. In \cite{bkr} a particular
crepant resolution is constructed as a moduli space.
It is expected that this is
a general phenomenon: every crepant resolution can be constructed as a moduli space.
This is proved in \cite{ci05} when $X =\CC^3/G$ where $G$ is abelian. We
speculate that the same is true
for global orbifolds: every crepant resolution of a global orbifold can be constructed
as a moduli space. Some variants of the Quot functor may be helpful. We hope to
return to this problem in  the future.
\end{Rem}
\begin{Rem}
Consider the case when $X$ has only symplectic quotient singularities. If $(\text{Hilb}'(\mX))^n \to X$ is a crepant resolution, then we have an equivalence of derived categories
 $D^b(\mX) \to D^b((\text{Hilb}'(\mX))^n)$,  cf \cite{kal05}.
\end{Rem}

In \cite{bp05}, a conjectural twisted version of derived McKay correspondence is
formulated. This conjecture suggests an equivalence between  the derived
category of twisted sheaves on an orbifold and a related derived category of
twisted sheaves on a crepant resolution. We establish this conjecture in the
setting of \cite{bkr}, see Theorem \ref{tbkr}.

\subsection{Differential graded categories}
It is known (e.g. \cite{bk90}) that the derived category $D(X)$ of coherent sheaves can be enriched to a differential graded (DG) category which we denote by $L_{coh}(X)$. Roughly speaking, objects in $L_{coh}(X)$ are complexes of sheaves. For two objects $A=(A^\cdot)$ and $B=(B^\cdot)$, a morphism $\phi$ of degree $i$ between them is a collection $(\phi^\cdot)$, $\phi^\cdot: A^\cdot \to B^{\cdot+i}$. The space $Hom_{DG}(A,B)$ of morphisms is a $\mathbb{Z}$-graded vector space endowed with a differential $d\phi=d_B\circ \phi-(-1)^i\phi\circ d_A$. In other words, $Hom_{DG}$ is a complex. One recovers $D(X)$ by localizing the homotopy category of $L_{coh}(X)$ with respect to the subcategory of acyclic complexes. For more details, see \cite{ke99}.

Here we add a note to the theory of D-equivalence that Fourier-Mukai type equivalence can be ``lifted'' to the level of DG-categories.

\begin{Lem}
Let $F: L_{coh}(X)\to L_{coh}(Y)$ be a DG functor which induces an equivalence of derived categories $D(X)\to D(Y)$. Then $F$ is a quasi-equivalence of DG categories.
\end{Lem}
\begin{proof}
Recall that a quasi-equivalence is a DG functor $F: L_{coh}(X)\to L_{coh}(Y)$ such that
\begin{enumerate}
\item
The morphism $Hom_{DG}(A,B)\to Hom_{DG}(FA,FB)$ is a quasi-equivalence of complexes;
\item
$F$ induces an equivalence $D(X)\to D(Y)$ of triangulated categories.
\end{enumerate}
Hence we only need to show the first condition, i.e. the cohomologies coincide. In degree zero, this is part of the derived equivalence: $$H^0Hom_{DG}(A,B)\simeq H^0Hom_{DG}(FA,FB).$$ On the other hand, we have $$H^iHom_{DG}(A,B)=H^0Hom_{DG}(A,B[i])$$ $$\simeq H^0Hom_{DG}(FA,FB[i])=H^iHom_{DG}(FA,FB).$$ This completes the proof.
\end{proof}
Since pushforward, pullback, and tensor product functors can be lifted to DG level, it is clear that a Fourier-Mukai functor comes from a DG functor $L_{coh}(X)\to L_{coh}(Y)$ which, according to the lemma above, is a quasi-equivalence. It's clear from the proof that in this context, having a quasi-equivalence at DG level yields no new information. However, the DG structure is sometimes necessary: for example, to construct the B-model potential of a Calabi-Yau manifold using a recent work of K. Costello \cite{co05}.

This paper is organized as follows. Several properties of the Quot functors we need are established in Section \ref{osquot}. Theorem \ref{gbkr} is proved in Section \ref{pfgbkr}. In Section \ref{twmckay} we prove some cases of a conjecture in \cite{bp05} concerning twisted derived McKay correspondence.
\section*{Acknowledgements}
We thank D. Abramovich, V. Baranovsky, D. Ben-Zvi, B. Keller, M. Olsson, and J. Starr, for many important discussions. Part of this work was done during the AMS summer institute in algebraic geometry in Seattle. We thank the organizers for hospitality and support.

\section{Quot functors after Olsson-Starr}\label{osquot}
In this section we discuss some properties of Quot functors for Deligne-Mumford
stacks following Olsson-Starr \cite{os}. We will focus on the case of quotients
of $\mO_\mX$. Let $Quot(\mO_\mX/\mX)$ (respectively $Quot(\mO_X/X)$) denote the
Quot functor associated to the sheaf $\mO_\mX$ over $\mX$ (respectively the
sheaf $\mO_X$ over $X$). According to \cite{Gr} and \cite{os}, these two
functors are representable by projective schemes which we denote by $\text{Quot}(\mO_\mX/\mX)$ and $\text{Quot}(\mO_X/X)$ respectively.

\begin{Def}
There is a morphism $\text{Quot}(\mO_\mX/\mX)\to \text{Quot}(\mO_X/X)$ defined as follows: The exact functor $\pi_*: Coh(\mX)\to Coh(X)$ yields a natural transformation
$Quot(\mO_\mX/\mX) \Rightarrow Quot(\mO_X/X)$ of Quot functors. Let $\text{Quot}(\mO_\mX/\mX)\to \text{Quot}(\mO_X/X)$ be the induced morphism between the corresponding schemes.
\end{Def}
Note that $X\subset \text{Quot}(\mO_X/X)$ is an irreducible component. Let
$\text{Hilb}'(\mX)\subset \text{Quot}(\mO_\mX/\mX)$ denote the irreducible
component containing the preimage of $X_{sm}$. Restricting the morphism $\text{Quot}(\mO_\mX/\mX)\to \text{Quot}(\mO_X/X)$ yields the morphism $\text{Hilb}'(\mX)\to X$.
Note that $\text{Hilb}'(\mX)$ is contained in the locus $\text{Quot}^1(\mO_\mX/\mX)$ which parametrizes quotients with Hilbert polynomial $1$.

\begin{Lem}\label{ghilb}
Assume that $\mX$ is of the form $[M/G]$ with $M$ quasi-projective and $G$ finite. Then $\text{Hilb}'(\mX)$ is isomorphic to the irreducible component of $\text{G-Hilb}(M)$ containing the free $G$-orbits.
\end{Lem}
\begin{proof}
There is a natural morphism $\text{Hilb}'(\mX)\to \text{G-Hilb}(M)$ defined as
follows: Given an object $\mZ\subset \mX\times S$ of $Hilb'(\mX)(S)$.
 The scheme
$$\mZ\times_{(\mX\times S)}(M\times S)\subset M\times S$$
 is naturally an $S$-family
of $G$-clusters: Clearly the group scheme $G\times S$ acts fiberwise on $\mZ\times_{\mX\times S}(M\times S)$. Also, since $M\times S\to \mX \times S$ is a principle $G$-bundle, the space $\Gamma(\mO_{\mZ\times_{\mX\times S}(M\times S)})$ coincides with the regular representation $\CC[G]\otimes \Gamma(\mO_S)$ of $G\times S$.
Thus this defines an object in $\text{G-Hilb}(M)(S)$.

It is easy to check that this morphism is a closed immersion, as one may recover $\mZ\subset \mX\times S$ by taking the stack quotient $[\mZ\times_{\mX\times S}(M\times S)/G]$. Moreover, let $\text{G-Hilb}(M)\to X=M/G$ be the morphism induced from the Hilbert-Chow morphism.
Then the following diagram commutes:
$$ \begin{array}{ccc}
 \text{Hilb}'(\mX)            & \to&  \text{G-Hilb}(M) \\
\downarrow     &  \swarrow  &  \\
 X.            & &  \end{array}$$
The result follows.
\end{proof}

An important property is that $\text{Hilb}'(\mX)$
behaves well under \'etale base-change
on the {\em coarse moduli space}.
\begin{Prop}\label{local}
Notation as above. Let $U\to X$ be
an \'etale morphism from a
scheme $U$. Then we have
$$U\times_X (\text{Hilb}'(\mX))^n\simeq (\text{Hilb}'(U\times_X\mX))^n.$$
\end{Prop}

\begin{proof}
To simply the notations, we  denote $\text{Hilb}'(\mX)$ by $W$,
$\text{Hilb}'(U\times_X\mX)$ by $W_U$. Let $W^n$ be the
normalization of $W$ and $W_U^n$ the normalization of $W_U$. 
Observe  that $W \times_X U$ and $W_U$ are birational (both of them
are birational to $U$).

Consider the universal family $\mZ \to W$ of quotient sheaves of $\mO_{\mX}$ (with Hilbert polynomial $1$), where 
$\mZ \subset W \times \mX$.
Note that $\mZ$ is indeed a substack of
$W \times_X \mX$. Pull back the universal 
family $\mZ \to W$ to the scheme $U \times_X W$. The resulting family is a
family of quotient sheaves of $\mO_{U \times_X \mX}$ with Hilbert
polynomial $1$.  Hence, there is a morphism $W \times_X U \to
W_U$. Note that $W \to X$ and $W_U \to U$ are  projective. The
base change of $W \to X$ to  $W \times_X U \to U$ is also
projective.
 Since 
$W$ is the fine moduli space (i.e. for $w_1 \neq w_2 \in W$, the corresponding quotient sheaves are distinct),    and $U \to X$ is quasi-finite, the
morphism $W \times_X U \to \text{Hilb}'(\mX \times_X U)=W_U$ is also
quasi-finite. Since both $W \times_X U$ and $W_U$ are projective
over $U$, it follows that $W \times_X U \to W_U$ is finite.

Since $W^n \to W$ is finite, it follows that $W^n \times_X U \to
W_U$ is finite. Since $W^n \times_X U \to W_U^n$ is finite and
birational, we conclude by Zariski's main theorem that $W^n
\times_X U \cong W_U^n$.
\end{proof}
Consider the universal quotient sheaf $\mZ \to \text{Hilb}(\mX)$ and the universal quotient sheaf $\mZ_{U} \to   
\text{Hilb}(\mX) \times_X U$. Pull  back  $\mZ$ to $(\text{Hilb}(\mX))^n$, denoted by $\mZ'$. Pull back $\mZ_{U}$ to 
$U \times_X (\text{Hilb}(\mX))^n$, denoted by $\mZ'_U$.
It can be seen from the proof that the further pull back of $\mZ'$ to  $\text{Hilb}(U\times_X \mX)$ 
is isomorphic to $\mZ'_U$.
Lemma \ref{ghilb} and Proposition \ref{local} imply in particular that in fact those G-Hilbert schemes  patch together nicely to a global object. 

\section{Proof of \ref{gbkr}}\label{pfgbkr}
In this section we present two proofs of Theorem \ref{gbkr}.
\subsection{Some basic materials}
We present several results concerning derived categories of coherent sheaves on smooth Deligne-Mumford stacks.

\begin{Prop}[Serre functor]\label{serre}
Let $\mX$ be a smooth separated Deligne-Mumford stack which has a coarse moduli space $X$ which is a quasi-projective Gorenstein
variety, whose dualizing sheaf is denoted by $\omega_X$.
Then $D^b(\mX)$ has a Serre functor $S(-):= (- \overset{\mathbf{L}}{\otimes} \pi^*\omega_X)[dim\mX]$.
\end{Prop}
\begin{proof}
Note that $S$ is clearly an equivalence. We need to show that for $u, v\in D^b(\mX)$ there is a bifunctorial isomorphism $$Hom(u,v)\simeq Hom(v, S(u))^\vee.$$ Our argument is parallel to that in \cite{ka02}, Proposition 2.6. First assume that $u$ is a locally free sheaf and $v$ is a sheaf with compact support. We may assume that $u=\mO_\mX$ by replacing $v$ by $u^\vee\otimes v$. Now we have $$Hom(u, v[k])\simeq H^k(\mX,v)\simeq H^k(X, \pi_*v).$$ On the other hand, $$Hom(v[k], S(u))\simeq Ext^{dim\mX-k}(v,\pi^*\omega_X)\simeq Ext^{dim \mX-k}(\pi_*v,\omega_X).$$ Hence $Hom(u, v[k])\simeq Hom(v[k],S(u))^\vee$ by the duality for $X$.

If $v$ is a locally free sheaf and $u$ is a sheaf with compact support, then by the previous case, we have
$$Hom(u, v[k])\simeq Hom(S(u),S(v)[k])\simeq Hom(v[k], S(u))^\vee.$$ The general case follows by taking locally free resolutions.
\end{proof}
Note that if $\mX$ is isomorphic to $X$ in codimension $1$, then $\pi^*\omega_X\simeq \omega_\mX$ as Cartier divisors. In this case the Serre functor is given by $S(-)=(- \overset{\mathbf{L}}{\otimes}\omega_\mX)[dim\mX]$.
\begin{Prop}[a spanning class]\label{spanning}
Let $\mX$ be a smooth Deligne-Mumford stack which has Serre duality. Then the set $$\{\mO_\mZ|\mZ\subset \mX \text{ is a closed substack, }\pi(\mZ) \text{ is a point in }X \}$$ is a spanning class of $D^b(\mX)$.
\end{Prop}
\begin{proof}
This follows from the argument of \cite{br99}, Example 2.2.
\end{proof}

It follows that a Deligne-Mumford stack as in Proposition \ref{serre} has a spanning class given as above.

\subsection{Reduction to local case: the first proof}
Let $\{U_i\}$ be an \'etale cover of $X$ such that $U_i\times_X\mX \simeq [M_i/G_i]$ for some schemes $M_i$ and finite groups $G_i$. Consider the diagram
$$\begin{CD}
 \mX \\
@V{\pi}VV  \\
X @< << (\text{Hilb}'(\mX))^n,
\end{CD}$$
and its pullback to $U_i$
$$\begin{CD}
\mX_i:= U_i\times_X\mX \\
@V{\pi_i}VV  \\
U_i @< << U_i\times_X(\text{Hilb}'(\mX))^n\simeq (U_i\times_X \text{Hilb}'(\mX))^n.
\end{CD}$$
Note that $(\text{Hilb}'(\mX) \times_X U_i)^n \cong (\text{Hilb}'(\mX \times_X U_i))^n$ (Proposition \ref{local}) and 
$\text{Hilb}'(\mX \times_X U_i)$ is already smooth.
Let $F: D^b((\text{Hilb}'(\mX))^n)\to D^b(\mX)$ be a Fourier-Mukai type transformation defined by an object $\mE\in D^b(\text{Hilb}'(\mX)\times_X \mX)$ and $F_i: D^b(U_i\times_X(\text{Hilb}'(\mX))^n)\to D^b(\mX_i)$ the Fourier-Mukai type transformation given by pulling everything back to $U_i$.

\begin{Prop}\label{glue}
Assume that $F_i$ are equivalences of derived categories for all $i$, then so is $F$.
\end{Prop}

\begin{proof}
This follows from the argument of \cite{ch02}, Proposition 3.2, using the spanning class given in Proposition \ref{spanning}.
\end{proof}

\begin{proof}[Proof of Theorem \ref{gbkr}]
By Lemma \ref{ghilb} and Proposition \ref{local}, the functor
$$F_i: D^b((U_i\times_X\text{Hilb}'(\mX))^n)\to D^b(\mX_i)$$ is the Fourier-Mukai type
transformation defined by the universal object of the Hilbert scheme
$\text{G}_i\text{-Hilb}(M_i)$. Again by Proposition \ref{local},  Lemma \ref{ghilb} and \cite{bkr}, $U_i\times_X(\text{Hilb}'(\mX))^n \cong 
(\text{Hilb}'(\mX \times_X U_i))^n \cong \text{Hilb}'(\mX \times_X U_i)$.
By the results of \cite{bkr}, we also know that $F_i$
is an equivalence, $\text{Hilb}'(\mX \times_X U_i)$ is smooth, and
$\text{Hilb}'(\mX \times_X U_i)\to U_i$ is a crepant resolution. It follows that $(\text{Hilb}'(\mX))^n$ is smooth and $(\text{Hilb}'(\mX))^n\to X$ is a
crepant resolution. By Proposition \ref{glue} it follows that $F$ is an equivalence.
\end{proof}

\subsection{Patching local crepant resolutions}
In general, without knowing the global space $\text{Hilb}'(\mX)$, it is not a priori clear the crepant resolutions of $U_i$ given by G-Hilbert schemes can be patched to a crepant resolution of $X$. We note that it in fact can be done in dimension 3.
\begin{Prop}\label{3dcr}
Let $X$ be a normal $\QQ$-factorial projective threefold and $\{U_i\}$ an open cover of $X$. Suppose that for every $i$ there is a crepant resolution $\phi_i: Y_i\to U_i$, then there is a crepant resolution $\phi: Y\to X$.
\end{Prop}
\begin{proof}
Consider a relative minimal model $f: Y\to X$, where $Y$ is terminal, $\QQ$-factorial and $K_Y$ is $f$-nef. We only need to prove that $Y\to X$ is crepant and $Y$ is smooth.
Note that $X$ has only canonical singularities since local crepant resolutions exist.
The first part follows since $X$ is canonical and $K_Y$ is $f$-nef. To check that $Y$ is smooth we only need to check it locally. Thus we consider the restriction $f|_{U_i}: Y|_{U_i}\to U_i$. Both $Y|_{U_i}\to U_i$ and $Y_i\to U_i$ are relative minimal models of $U_i$, so they can be connected by a sequence of flops. Since flops in dimension 3 preserve smoothness, the result follows.
\end{proof}

\begin{Rem}
This argument does not work in higher dimensions, since flops may not terminate and may
not preserve smoothness.
\end{Rem}

\subsection{Another proof of \ref{gbkr}}\label{another}
Instead of using results of \cite{bkr}, the second proof uses their arguments. It is not hard to check that the arguments in \cite{bkr} extends to our case provided Proposition \ref{local} and the following:
\begin{enumerate}
\item
the category $D^b(\mX)$ is indecomposable;
\item
we have the Grothedieck duality for the morphism $(\text{Hilb}'(\mX))^n\times \mX \to (\text{Hilb}'(\mX))^n$;
\item
we have a Serre functor for $D^b(\mX)$.
\end{enumerate}

(1) follows from the argument of \cite{br99}, Example 3.2, provided we know that for any integral closed substack $\mW$ of $\mX$, the sheaf $\mO_\mW$ is indecomposible. This can be seen as follows: Since $\mX$ is smooth, we have $\mX\simeq [M/GL_r]$ where $GL_r$ is some general linear group. By assumption $\mX$ is connected, so is $X$. It follows that $M$ is connected. Let $p:M \to \mX$ be the structure morphism. Then $p^*: D^b(\mX)\to D^b(M)$ is fully faithful, and $p^*\mO_\mW=\mO_{\mW\times_\mX M}$ is indecomposible. Hence $\mO_\mW$ is indecomposible.

For (2), note that $(\text{Hilb}'(\mX))^n\times \mX \to (\text{Hilb}'(\mX))^n$ factors as $$(\text{Hilb}'(\mX))^n\times \mX \overset{id\times \pi}{\longrightarrow} (\text{Hilb}'(\mX))^n\times X \to(\text{Hilb}'(\mX))^n.$$ Also note that $(id\times \pi)_*$ is exact and has the left and right adjoint $(id\times \pi)^*$. We may conclude by the Grothendick duality for schemes.

(3) follows from Proposition \ref{serre}.

\section{Twisted derived McKay correspondence}\label{twmckay}
In this section we discuss the twisted version of derived McKay correspondence proposed in \cite{bp05}. Let $X=M/G$ with $M$ quasi-projective of dimension $n$ and $G$ finite and $\mX$ be the
associated stack as discussed above. Also let $Y$ be the component of $G$-Hilbert scheme of $M$ which contains free orbits. According to \cite{bkr}, we have an integral functor
$\Phi: D^b(Y)\to D^b(\mX)$ given by the kernel $\mO_{\mZ}$, which is an equivalence under certain conditions. Here $\mZ=[Z/G]$ and $Z\subset Y\times M$ is the universal closed subscheme.

The authors of \cite{bp05} proposed a twisted version of this equivalence, which we recast in our setting as follows: Let $Br(Y)$ be the Brauer group of $Y$. It is shown in
\cite{bp05} that both $Br(Y)$ and $H^2(G, \CC^*)$ can be embedded into a larger group, namely $Br(X_{sm})$. Let $\alpha \in Br(Y)\cap H^2(G,\CC^*)$ be a class which is
$r$-torsion. It follows that there is an injective homomorphism $\iota: \mu_r\to\CC^*$ such that $\alpha$ is in the images of the induced maps $H_{et}^2(Y, \mu_r)\to H_{et}^2(Y,\CC^*)$ and $H^2(G,\mu_r)\to H^2(G,\CC^*)$. Thus we can associate a $\mu_r$-gerbe $\mY_\alpha$ over $Y$ and a $\mu_r$-gerbe $\mX_\alpha$ over $\mX$. The $\alpha$-twisted derived categories $D^b(Y, \alpha)$ and
$D^b(\mX, \alpha)$ are derived categories of coherent sheaves on $\mY_\alpha$ and $\mX_\alpha$ whose actions by $\CC^*$ and $\mu_r$ are compatible via $\iota: \mu_r\to\CC^*$ (see
\cite{li04}). It is conjectured in \cite{bp05} that there is a equivalence of derived categories $$D^b(Y, \alpha)\to D^b(\mX,\alpha).$$
We prove this conjecture in our situation. Consider the following diagram:
$$\begin{CD}
\mX_\alpha\times\mY_\alpha@>\pi_Y>>\mY_\alpha@>p_Y>>Y\\
@V{\pi_\mX}VV @. @V \bar{\phi}VV\\
\mX_\alpha@>p_\mX>>\mX@>p_X>> X
\end{CD}$$

Since $\alpha\in Br(Y)\cap H^2(G,\CC^*)$, it is not hard to check that the integral functor $\Phi_\alpha$ defined by the kernel $(p_Y\times p_\mX)^*\mO_{\mZ}$ is a functor between
$D^b(Y,\alpha)$ and $D^b(\mX,\alpha)$.
\begin{Thm}\label{tbkr}
Assume that $dim Y\times_X Y\leq n+1$, then $\Phi_\alpha$ is an equivalence of triangulated categories.
\end{Thm}
\begin{proof}
According to \cite{bkr}, $Y$ is smooth and $Y\to X$ is a crepant resolution. It
follows that $\mY_\alpha$ and $\mX_\alpha$ are smooth Deligne-Mumford stacks.
In view of this, that $\Phi_\alpha$ is an equivalence is a special case of a
stack version of Bridgeland's result (\cite{br99}, Theorem 1.1). It is
straightforward (although lengthy) to modify Bridgeland's arguments in
\cite{br99} to prove that $\Phi_\alpha$ is an equivalence if and only if
\begin{equation}\label{breq1}
Ext_{\mX_\alpha}^i((p_Y\times p_\mX)^*\mO_\mZ|_{y_1}, (p_Y\times p_\mX)^*\mO_\mZ|_{y_2})=0 \text{ for all } i\geq 0, y_1\neq y_2\in \mY_\alpha;
\end{equation}

\begin{equation}\label{breq2}
(p_Y\times p_\mX)^*\mO_\mZ|_{y}\simeq (p_Y\times p_\mX)^*\mO_\mZ|_{y}\otimes p_\mX^*p_X^*\omega_X \text{ for all }y\in \mY_\alpha.
\end{equation}

We simply note that in order for the arguments in \cite{br99} to work for this, we
need to have a Serre functor, a spanning class and indecomposability of the
derived categories, as well as calculations of certain $Ext$ groups. These have
been settled in Proposition \ref{serre}, \ref{spanning} and in Section
\ref{another}.
The needed calculations of $Ext$ groups follows from those in
\cite{bkr}. Now the calculations in \cite{bkr} immediately imply (\ref{breq1})
and (\ref{breq2}). Hence $\Phi_\alpha$ is an equivalence.
\end{proof}

\begin{Rem}
\hfill
\begin{enumerate}
\item
It is natural to ask for a global version of the conjecture in \cite{bp05}. One expects
that the work of \cite{bp05} can be generalized to give such a conjecture. Our
work should be helpful in proving such a conjecture on global twisted derived McKay
correspondence.

\item
Consider the following situation: $X=M/G$ as above,  $\mX=[M/G]$ is the Deligne-Mumford stack as before, and $Y\to X$ be a crepant resolution. Suppose that there is an equivalence of derived categories $D^b(Y)\simeq D^b(\mX)$. It is known that such an equivalence is a Fourier-Mukai transform given by an integral kernel $\mE\in D^b(Y\times \mX)$. If $\mE$ is a sheaf on $Y\times \mX$, then our arguments can be modified to show that for $\alpha\in Br(Y)\cap H^2(G, \CC^*)$, the integral functor given by the pullback of $\mE$ to the associated gerbe yields an equivalence $D^b(Y,\alpha)\simeq D^b(\mX,\alpha)$. It is interesting to understand whether this is the case or not, if $\mE$ is not a sheaf.
\end{enumerate}
\end{Rem}

\end{document}